\documentclass[12pt]{amsart}
\usepackage{graphicx}
\usepackage{amsmath,amssymb,amsthm}
\usepackage[final]{showkeys}  
\usepackage{mathbbol,bbm}
\usepackage[english]{babel}

\newcommand{\Rl}{\mathbb{R}}

\newcommand{\Aff}{\mathbb{A}}
\newcommand{\Euc}{\mathbb{E}}

\newcommand{\pard}[2]{\frac{\partial #1}{\partial #2}}

\DeclareMathOperator{\grad}{grad}

\newcommand{\scal}[1]{\langle #1\rangle}
\newcommand{\lie}[1]{\left[#1\right]}
\newcommand{\lc}[1]{\hat{\nabla}_{#1}}
\providecommand{\abs}[1]{\lvert#1\rvert}
\providecommand{\norm}[1]{\lVert#1\rVert}

\newtheorem{lemma}{Lemma}[section]
\newtheorem{theorem}{Theorem}[section]

\begin{document}

\title{Affine hypersurfaces admitting a pointwise $SO(n-1)$ symmetry.}

\author[K. Schoels]{Kristof Schoels}
\email[K. Schoels]{kristof.schoels@wis.kuleuven.be}
\address[K.Schoels]{Department of Mathematics\\ Katholieke Universiteit Leuven\\
Celestijnenlaan 200 B, Box 2400\\BE-3001 Leuven\\ Belgium}

\begin{abstract}
In this article we obtain a classification of strictly locally convex affine hypersurfaces in $\Aff^{n+1}$ for which the geometrical structure is pointwise invariant under the group $SO(n-1)$ represented by rotations around a fixed axis in the tangent space. This generalises the results obtained by Scharlach and Lu for $n=3$ in \cite{scharl}.
\end{abstract}
\maketitle

\section{Introduction.}
An immersion $\phi:M\rightarrow N$ of a submanifold $M$ in an ambient space $N$ defines at each $p\in M$ a linear monomorphism $\phi_*:T_pM\rightarrow T_{\phi(p)}N$. This identifies the tangent space to $M$ with a subspace of the tangent space to $N$. One can look for a suitable subspace $V_p \subset T_{\phi(p)}N$ which is complementary to this subspace. Requiring differentiability, one aims to construct a transversal distribution to complement the tangential distribution. When $N$ is a Riemannian manifold with a metric $g$, such an immersion $\phi$ induces a Riemannian structure on $M$. It is common to choose the orthogonal distribution in respect to $g$ as the transversal distribution. The geometrical structure resulting through this construction is invariant under isometries of $N$. To allow for invariance under a wider range of transformations, other normalisations are in order.\\
When we consider hypersurfaces $M\subset\Aff^{n+1}$, we can take any transversal vector field $\xi$ and decompose, for $X,Y$ tangent to $M$, the natural covariant derivative into a transversal and tangential part yielding
\begin{align*}
D_XY&=\nabla_XY+h(X,Y) \xi,\\
D_X \xi&= -SX+\tau(X) \xi.
\end{align*}
The induced structure $(S,\nabla,h,\tau)$ determines the immersion completely up to an affine transformation of $\Aff^{n+1}$. This structure does depend on the choice of transversal vector field $\xi$, but under the transformation $\xi\rightarrow\lambda \xi+Z$, the bilinear symmetric form $h$ changes to $\frac{h}{\lambda}$, hence degeneracy and signature are independent of $\xi$. Another object that is of importance is volume. The standard volume form in $\Aff^{n+1}$ obtained by taking the determinant is invariant under equi-affine transformations. For tangential vector fields $\{X_1,\cdots,X_n\}$ the induced volume form $\omega$ on $M$, dependent on $\xi$, is given by
\begin{equation*}
\omega(X_1,X_2,\cdots,X_n)=\begin{vmatrix} \phi_*(X_1) &\cdots & \phi_*(X_n) & \xi \end{vmatrix}.
\end{equation*}
Under the previous change of transversal $\xi$, this volume form changes as $\omega\rightarrow \lambda \omega$. The following result gives a natural choice for a transversal vector field on non-degenerate hypersurfaces, see \cite{nomsaz}.
\begin{theorem}[Blaschke]
Suppose $\phi:M^n\rightarrow \Aff^{n+1}$ is an immersion of a connected hypersurface, for which the induced bilinear form $h$ is non-degenerate. Then there exists a transversal vector field $\xi$ such that
\begin{enumerate}
\item $\nabla \omega=0$,\\
\item $\omega(X_1,\cdots,X_n)^2=\abs{\det(h_{ij})}$.
\end{enumerate}
This vector field is unique up to sign. This vector field is referred to as the affine normal.
\end{theorem}
The first demand is equivalent to demanding that $\tau=0$. Because $h$ is non-degenerate, $\omega$ defines a volume form. When $h$ is definite, the ambiguity concerning the sign can be removed by requiring that $h$ is positive definite. In this case, the manifold curves towards $\xi$ and is strictly locally convex. Graphically, the affine normal in a point $p\in M$ can be obtained by slicing $M$ with hyperplanes parallel to $T_pM$. The barycenters of the resultant $(n-1)$-manifolds form a curve $C:I\subset\Rl\rightarrow \Aff^{n+1}$. The tangent direction to this curve in $p$ coincides with the direction of the affine normal, see \cite{simonlizhao}.\\
Using the affine normal, the manifold has a Blaschke connection $\nabla$, an affine metric $h$ and a shape operator $S$. Calculating the curvature tensor of the flat manifold $\Aff^{n+1}$, the following equations hold:
\begin{align*}
\nabla h\quad &\mbox{is symmetric},\\
\nabla S\quad &\mbox{is symmetric},\\
h(X,SY)&=h(Y,SX),\\
R(X,Y)Z&=h(Y,Z)SX-h(X,Z)SY.
\end{align*}
If $n\geq 3$, the symmetry of $\nabla S$ and $\nabla h$ can be shown to be equivalent, using the second Bianchi relation for $R$, see \cite{nomsaz}. These equations are sufficient for a manifold $M$ with structure $(S,\nabla,h)$ to be immersed as an affine hypersurface, uniquely up to an affine transformation. A particular case of affine hypersurfaces arises when $S=\lambda \Eins$. These hypersurfaces are referred to as affine spheres. These are the hypersurfaces for which the affine normals are parallel in $\Aff^{n+1}$ (when $S=0$) or the affine normals intersect in a common center (when $S\neq0$). For Euclidean hypersurfaces this property implies that the hypersurface is either part of a Euclidian sphere or a part of a plane. However, a similar result no longer holds for affine hyperspheres. Examples of these spheres are quadrics, such as the ellipsoid, the one sheeted hyperboloid and the elliptic paraboloid. Another example is the Calabi sphere given implicitly by $\prod_{i=1}^{n+1}x_i=1$. The Calabi conjecture shows that even under a global condition of completeness, every convex hypercone has a unique homothetic family of hyperbolic affine spheres asymptotic to its boundary, see \cite{simonlizhao}. So the amount of non-equivalent affine hyperspheres is copious.\\
Since $h$ is a metric, $(M,h)$ can be considered a (semi)-Riemannian manifold equipped with the Levi-Civita connection $\hat{\nabla}$. Comparing $\nabla$ and $\hat{\nabla}$ yields a symmetric $(1,2)$-tensor $K$ given by
\begin{equation*}
K(X,Y)=\nabla_XY-\hat{\nabla}_XY.
\end{equation*}
Using the symmetry of $\nabla h$, this difference tensor satisfies
$$h(K(X,Y),Z)=h(K(X,Z),Y)$$
and since $\hat{\nabla} \omega=\nabla \omega=0$, the tensor $K_X$ is traceless for every choice of $X$. One has the following theorem, see \cite{nomsaz}.
\begin{theorem}[Pick-Berwald]
An affine hypersurface $M$ for which $K=0$ is a part of a quadratic hypersurface.
\end{theorem}
For convex hypersurfaces, these are the ellipsoid, paraboloid and hyperboloid. We will refer to them as $\epsilon$-quadrics, where $\epsilon$ is the mean curvature. Thus $K$ plays an important part in describing the submanifold. Using it to connect the curvature tensor $\hat{R}$ of $\hat{\nabla}$ to $R$ of $\nabla$, we obtain the Gauss equation
\begin{equation}
\begin{split}
\hat{R}(X,Y)Z=\frac{1}{2}&\big(h(Y,Z)SX-h(X,Z)SY+h(SY,Z)X-h(SX,Z)Y\big)\\
            &\qquad-\lie{K_X,K_Y}Z\label{gauss}
\end{split}
\end{equation}
and the Codazzi equation
\begin{equation}
\begin{split}
\hat{\nabla}K(X;Y,Z)-\hat{\nabla}K(Y;X,Z)&=\frac{1}{2}\big(h(Y,Z)SX-h(X,Z)SY\\
&\qquad+h(SX,Z)Y-h(SY,Z)X\big).\label{codazzi}
\end{split}
\end{equation}
When we add the Ricci equation
\begin{equation}
\lc{} S(X;Y)-\lc{} S(Y;X)=K(SX,Y)-K(SY,X), \label{ricci}
\end{equation}
we have enough equations to determine an affine hypersurface completely, see \cite{scharl}.
\begin{theorem}
Suppose $M^n$ is a manifold equipped with a metric $h$. Suppose $S$ is an $h$-symmetric $(1,1)$-tensor and $K$ is a traceless $h$-symmetric $(1,2)$-tensor. Suppose that the Levi-Civita connection $\hat{\nabla}$, $K$ and $S$ satisfy equations \eqref{gauss} to \eqref{ricci}. Then there exists an immersion
\begin{equation*}
\phi:M\rightarrow \Aff^{n+1}
\end{equation*}
as an affine hypersurface with affine metric $h$, shape operator $S$ and difference tensor $K$. This immersion is unique up to equi-affine congruence.
\end{theorem}
 The question we try to solve is, given $K$ and $S$ what do the equations reveal about $\hat{\nabla}$ and in what way can the manifold be immersed as an affine hypersurface?\\
It is a tedious task to look for affine hypersurfaces with arbitrary $K$ and $S$. Following the idea of Bryant in \cite{bryant}, we impose certain pointwise symmetries to reduce the amount of variables we deal with. In this case, we assume that $K$ and $S$ are pointwise invariant under rotations around a given tangent vector $X_1$. We will show that this implies that for a well chosen $h$-orthonormal frame $\{X_1,\cdots,X_n\}$, $S$ and $K$ are given as
\begin{equation}
\begin{split}
S&=\begin{pmatrix} a & 0\\0 & b\Eins_{n-1}\end{pmatrix},\\
K_{X_1}&=\begin{pmatrix} (n-1)r & 0\\0 & -r \Eins_{n-1} \end{pmatrix},\\
(K_{X_i})_{AB}&=-r\left(\delta_{Ai}\delta_{B1}+\delta_{Bi}\delta_{A1}\right).
\end{split}\label{son-1}
\end{equation}
We will use the indices $\{A,B,C\}$ to run from $1$ to $n$ and $\{i,j,k,l\}$ when we leave out $\{1\}$. We will see that the entries $a$, $b$ and $r$, initially defined pointwise, are in fact smooth functions along the hypersurface and $r$ vanishes nowhere. We will see that the resulting hypersurface is a warped product $\Rl\times_{e^f} N^{n-1}$, where $N$ is a quadric. The warping occurs through a planar equi-affine curve subject to appropriate conditions.\\
First we elaborate on the immediate implications of the $SO(n-1)$ symmetry. Then we present the full classification. We conclude by specifying the results with the added condition that the hypersurface is an affine sphere and we give a few examples. The results are built upon classifications of $SO(2)$ symmetry in 3-dimensional affine hypersurfaces and spheres, see \cite{scharl} and \cite{vranck}. Certain special cases of these hypersurfaces can be constructed as Calabi-type compositions found in \cite{dill}.

\section{Preliminaries}
Before looking for the actual immersions, we show how the $SO(n-1)$ group shapes the tensors $S$ and $K$ and what this implies for \eqref{gauss} to \eqref{ricci}.\\
A pointwise local symmetry group on $M$ assigns for every $p\in M$ a subgroup $G$ of $Aut(T_pM)$ such that for every choice $g\in G$ one has
\begin{align*}
\forall X,Y \in T_pM&:h(gX,gY)=h(X,Y),\\
\forall X\in T_pM&: S(gX)=gS(X),\\
\forall X,Y \in T_pM&:K(gX,gY)=gK(X,Y).
\end{align*}
We assume that for every point $p\in M$ the symmetry is the same. Notice that for $h$ a definite metric, $G$ has to be a subgroup of the orthogonal group $O(n)$. There are a plethora of possible subgroups of $O(n)$, especially as $n$ increases. We consider one special case, the subgroup $SO(n-1)$, which manifests itself as the group of rotations around a fixed axis. Let us therefore consider a point $p$ and a unit tangent vector $X_1$, which represents the rotational axis. It splits the tangent space in $2$ parts: $T_pM=\scal{X_1}\oplus \scal{X_1}^{\bot}$. The group $SO(n-1)$ is given by
\begin{equation*}
SO(n-1)=\{g\in Aut(T_pM)\| gX_1=X_1, g:\scal{X_1}^\bot\rightarrow \scal{X_1}^\bot \in SO(n-1,\scal{X_1}^\bot)\}.
\end{equation*}
As a representation, it splits into a trivial representation on a $1$-dimensional space and the irreducible standard representation on an $(n-1)$-dimensional space. We can find the following result.
\begin{lemma}
Suppose $T_pM$ is a vector space of dimension $n\geq 3$ and $h$ is a positive definite metric on $T_pM$, $S$ is a $(1,1)$-tensor and $K$ is a symmetric $(1,2)$-tensor. Suppose that both $S$ and $K$ are $h$-symmetric and $K$ is traceless. Lastly, suppose that $S$ and $K$ are invariant under rotations about the axis determined by the unit normal $X_1$. Then for an orthonormal basis $\{X_1,X_2,X_3,\cdots,X_n\}$, the maps $S$ and $K$ are given by \eqref{son-1} for $r,a,b\in \Rl$.
\end{lemma}
\begin{proof}
First, let's consider the map $S$. We can write $SX_1=aX_1+V$ where $V\in\scal{X_1}^\bot$. Applying an element $g\in SO(n-1)$ before and after $S$, we conclude
\begin{align*}
SgX_1=SX_1=aX_1+V=aX_1+OV=gSX_1.
\end{align*}
Hence $V$ is invariant under the rotation $O$. Because this is true for any rotation $O$, $V=0$. Using that $S$ is $h$-symmetric, the component of $SX$ along $X_1$ vanishes when $h(X,X_1)$ does. So $S$ restricts to an endomorphism of $\scal{X_1}^\bot$. This endomorphism commutes with the irreducible action of $SO(n-1)$ on this subspace, hence applying Schur's lemma we conclude that this reduced map should be a multiple of the identity. For details on this, see \cite{schur}. Hence, we find a number $b$ such that
$$S=\begin{pmatrix} a & 0\\0 & b\Eins_{n-1}\end{pmatrix}.$$
The same argument shapes the operator $K_{X_1}$. We just have to look at the shape of $K(X,Y)$ for $X,Y\in \scal{X_1}^\bot$. We can write this as $K(X,Y)=-r h(X,Y) X_1+K_2(X,Y)$ for which $K_2$ satisfies the same conditions as $K$, but on the space $\scal{X_1}^\bot$ and it is invariant under every rotation. We have to proof that $K_2=0$. Because $K_2$ is symmetric, it is sufficient to proof that $K_2(X,X)=0$ for every choice of $X\in \scal{X_1}^\bot$. Suppose this is not true for a certain $X$. Choose an orthogonal transformation $O$ such that $OX=-X$. Hence we get that $OK_2(X,X)=K_2(X,X)$. Now specify $O$ such that $OK_2(X,X)=-K_2(X,X)$. This can be achieved by splitting $K_2(X,X)$ into a component in the direction of $X$ and $Y$ orthogonal to $X$. Hence, we specify $OY=-Y$. We can always complete $\{X,Y\}$ to an orthogonal basis of $T_pM$ and let $O$ map each other vector in the basis to itself, so $O$ is an element of $SO(n-1)$. If $Y=0$, it is sufficient to flip the sign of just one vector orthogonal to $X$ in order to obtain an element of $SO(n-1)$. Hence $K_2(X,X)=-K_2(X,X)=0$, which contradicts our initial assumption. This gives $K$ the desired form.
\end{proof}
We can assume that $r> 0$, at least on a dense open subset. If $r=0$ on an open neighbourhood of $M$, this neighboorhood would be a part of a quadric and so the local symmetry group would be $O(n)$. This leaves us with vectors $\{X_1,X_2,\cdots,X_n\}$ and numbers $a$, $b$ and $r$ defined pointwise. We would like to have differentiable functions and vector fields overall. To obtain this, we will assume that \eqref{son-1} is satisfied pointwise on the manifold. Using \eqref{gauss} we can calculate the curvature tensor $\hat{R}$, which is defined pointwise. We obtain that
\begin{align*}
\hat{R}(X_1,X_i)X_1&=-\left(nr^2+\frac{a+b}{2}\right)X_i,\\
\hat{R}(X_1,X_i)X_j&=\delta_{ij}\left(\frac{a+b}{2}+nr^2\right)X_1,\\
\hat{R}(X_i,X_j)X_1&=0,\\
\hat{R}(X_i,X_j)X_k&=(b-r^2)\left(\delta_{jk}X_i-\delta_{ik}X_j\right).
\end{align*}
Its Ricci tensor can be constructed and results in
\begin{equation*}
\hat{Ric}=\begin{pmatrix} (n-1)\left(nr^2+\frac{a+b}{2}\right) & 0\\
0 & \left(\frac{a}{2}+\frac{2n-3}{2}b+2r^2\right)\Eins_{n-1}\end{pmatrix}
\end{equation*}
The Riccitensor has $2$ Eigenvalues. Those coincide when
\begin{equation*}
(n-1)\left(nr^2+\frac{a+b}{2}\right)=\left(\frac{a}{2}+\frac{2n-3}{2}b+2r^2\right)\Leftrightarrow \frac{a-b}{2}=-(n+1)r^2.
\end{equation*}
If $a=b$, this occurs when $r=0$ which we have assumed not to be the case. Otherwise, $X_1$ is the unique principal direction of $S$ with Eigenvalue equal to $a$, hence it is a differentiable vector field. When the Eigenvalues of $\hat{Ric}$ are distinct, the direction of the first Eigenvalue defines a differentiable vector field. When $X_1$ is the desired differentiable vector field, the orthogonal complement forms a differentiable distribution. Because the vector fields $\{X_1,\cdots,X_n\}$ are differentiable, the pointwise defined numbers $a$, $b$ and $r$ turn into differentiable functions too.\\
Now we have obtained an orthonormal frame $\{X_1,\cdots,X_n\}$ where $X_1$ determines the axis of rotation and $S$ and $K$ are determined by differentiable functions $a$, $b$ and $r>0$. We can apply \eqref{gauss} to \eqref{ricci} and obtain the following result.
\begin{lemma}
Suppose that $M$ is an affine hypersurface with a pointwise $SO(n-1)$ symmetry. Then for a natural frame $\{X_1,\cdots,X_n\}$ satisfying \eqref{son-1}, we obtain that
\begin{equation}
\begin{split}
\lc{X_1}X_1&=0,\\
\lc{X_i}X_1&=\sigma X_i,\\
X_i(a)&=X_i(b)=X_i(\sigma)=X_i(r)=0,\\
X_1(b)&=(a-b)\left(\sigma-r\right),\\
X_1(r)&=-\left(\frac{a-b}{2}+(n+1)\sigma r\right),\\
X_1(\sigma)&=-\left(\frac{a+b}{2}+nr^2+\sigma^2\right).\label{diffeq}
\end{split}
\end{equation}
\end{lemma}
\begin{proof}
We denote the components of the connection $\hat{\nabla}$ as
\begin{equation*}
\lc{X_A}X_B=\sum_{C=1}^n \Gamma_{AB}^C X_C.
\end{equation*}
We start by calculating \eqref{ricci} and find
\begin{align*}
\hat{\nabla} S(X_1;X_i)-\hat{\nabla} S(X_i;X_1)&=\left((a-b)\Gamma_{11}^i-X_i(a)\right) X_1+X_1(b) X_i-(a-b)\sum_{j=2}^n \Gamma_{i1}^jX_j, \\
K(SX_1,X_i)-K(SX_i,X_1)&=-r(a-b)X_i,\\
\hat{\nabla} S(X_i;X_j)-\hat{\nabla} S(X_j;X_i)&=(a-b)\left(\Gamma_{i1}^j-\Gamma_{j1}^i\right) X_1+X_i(b) X_j-X_j(b) X_i,\\
K(SX_i,X_j)-K(SX_j,X_i)&=0.
\end{align*}
This results in
\begin{align*}
X_i(a)&=(a-b) \Gamma_{11}^i,\\
X_i(b)&=0,\\
X_1(b)&=(a-b) (\Gamma_{i1}^i-r),\\
(a-b)\Gamma_{i1}^j&=0 \qquad \mbox{when $i\neq j$}.
\end{align*}
These equations are helpful, as long as the affine hypersurface is not a hypersphere. Otherwise, equation \eqref{ricci} implies that $S=a\Eins$ where $a$ is a constant. Secondly, we calculate \eqref{codazzi} and this shows
\begin{align*}
\hat{\nabla}K(X_1;X_i,X_i)&=-X_1(r)X_1-2r\Gamma_{11}^i X_i-r\sum_{j=2}^n\Gamma_{11}^jX_j,\\
\hat{\nabla}K(X_i;X_1,X_i)&=(n+1)r\Gamma_{i1}^iX_1-X_i(r)X_i,\\
\hat{\nabla}K(X_i;X_1,X_1)&=(n-1)X_i(r)X_1+(n+1)r\sum_{j=2}^n\Gamma_{i1}^jX_j,\\
\hat{\nabla}K(X_1;X_1,X_i)&=(n+1)r\Gamma_{11}^i X_1-X_1(r) X_i,\\
\hat{\nabla}K(X_1;X_i,X_j)&=-r\left(\Gamma_{11}^i X_j+\Gamma_{11}^j X_i\right),\\
\hat{\nabla}K(X_i;X_1,X_j)&=(n+1)r\Gamma_{i1}^j X_1-X_i(r) X_j,\\
\hat{\nabla}K(X_i;X_j,X_j)&=-X_i(r)X_1-2r\Gamma_{i1}^j X_j-r\sum_{k=2}^n \Gamma_{i1}^kX_k,\\
\hat{\nabla}K(X_j;X_i,X_j)&=-r\left(\Gamma_{j1}^iX_j+\Gamma_{j1}^j X_i\right),\\
\hat{\nabla}K(X_i;X_j,X_k)&=-r\left(\Gamma_{i1}^j X_k+\Gamma_{i1}^k X_j\right).
\end{align*}
Filling in the right-hand side of \eqref{codazzi} and comparing yields
\begin{align*}
\lc{X_1}X_1&=0,\\
\lc{X_i}X_1&=\Gamma_{i1}^iX_i,\\
X_i(r)&=0,\\
X_1(r)&=-\left(\frac{a-b}{2}+(n+1)r \Gamma_{i1}^i\right).
\end{align*}
Since $\Gamma_{i1}^i$ is the same for every choice of $i$, we relabel it $\sigma$. Finally applying \eqref{gauss} we obtain
\begin{align*}
h\left(\hat{R}(X_1,X_i)X_1,X_i\right)&=X_1(\sigma)+\sigma^2\\
                                     &=-\left(\frac{a+b}{2}+nr^2\right),\\
h\left(\hat{R}(X_i,X_j)X_1,X_j\right)&=X_i(\sigma)\\
&=0.
\end{align*}
This completes the proof.
\end{proof}
We can remark that there are other relations hidden inside the Gauss equations. But we will see that they are not necessary to construct all the immersions.\\
Before classifying the hypersurfaces, we give a brief description of warped products. Consider Riemannian manifolds $(M_1,g_1)$ and $(M_2,g_2)$. We can combine both manifolds as the set $M_1\times M_2$ of couples of points in $M_1$ and $M_2$. On this product manifold we define a Riemannian metric $g=g_1+e^{2f} g_2$, where $f$ is a function on $M_1\times M_2$. Given tangent vectors $X$ and $Y$, we can decompose them into a part tangent to $M_1$ and one tangent to $M_2$, using the projection maps to the respective manifold. Then $g(X,Y)$ can be calculated as
\begin{equation*}
g(X,Y)=g_1(X_1,Y_1)+e^{2f}g_2(X_2,Y_2).
\end{equation*}
Because $g_1$ and $g_2$ are Riemannian and $e^{2f}>0$, the resultant manifold $(M_1\times M_2,g)$ is a Riemannian manifold, referred to as the twisted product of $M_1$ and $M_2$ with twisting function $e^f$. In case $f$ only depends on $M_1$, the product is referred to as a warped product. An example of a warped product surface in $\Euc^3$ is a rotational surface of the planar curve $(\gamma_1,\gamma_2)$ around the $x$-axis. This produces the warped product $\Rl\times_{\gamma_2} S^1$. For twisted and warped products one has the following result, see \cite{hiepko} and \cite{ponge}.
\begin{theorem}
Given a Riemannian manifold $(M,g)$ with $2$ orthogonal distributions $\mathcal{N}_1$ and $\mathcal{N}_2$. Suppose that for any vector fields $X,Y\in \mathcal{N}_1$ and $U,V\in \mathcal{N}_2$ one has
\begin{align*}
\lc{X}Y &\in \mathcal{N}_1,\\
\lc{U}V \mod \mathcal{N}_2&=g(U,V) H,
\end{align*}
then $M$ is obtained locally as a twisted product $M_1\times_{e^{f}} M_2$, where $M_i$ is the integral manifold of $\mathcal{N}_i$ through the given point $p$. If $H=\norm{H} H_0$ and $U(\norm{H})=0$ for every $U\in \mathcal{N}_2$, then $M$ is locally a warped product $M_1\times_{e^f} M_2$ and $f$ is given as
\begin{equation*}
\grad{f}=-H.
\end{equation*}
\end{theorem}
As a corollary the hypersurfaces we consider are warped products $\Rl\times_{e^f} N^{n-1}$ and $X_1(f)=\sigma$. Furthermore, if $g_N$ is the metric on $N$, then we obtain
$$g_N(X_i,X_j)=e^{-2f}\delta_{ij}.$$
To obtain an orthonormal tangent frame on $N$, we can take $U_i=e^f X_i$. We can use the relationship between curvature tensors on a warped product from \cite{oneill} to calculate the sectional curvature of $(N,g_N)$ along a plane through $\{X_i,X_j\}$ and obtain
\begin{align*}
\kappa_N(U_i,U_j)&=e^{2f} h\left(\hat{R}(X_i,X_j)X_j,X_i\right)-e^{2f}\frac{\norm{\grad{e^f}}^2}{e^{2f}} \left(h(X_i,X_j)^2-\norm{X_i}^2\norm{X_j}^2\right)\\
                 &=e^{2f} (b-r^2+\sigma^2).
\end{align*}
This means that $N$ is a manifold of constant sectional curvature. As $\kappa_N$ can still change along $X_1$, we can calculate its derivative and we find
\begin{equation*}
X_1(\kappa_N)=2X_1(f) \kappa_N+e^{2f}\left(X_1(b)-2rX_1(r)+2\sigma X_1(\sigma)\right)=2\sigma\kappa_N-2\sigma\kappa_N=0.
\end{equation*}
Hence we can rescale $N$ to obtain $\kappa_N=\epsilon\in\{-1,0,1\}$. Notice that $\zeta=b-r^2+\sigma^2=e^{-2f}\epsilon$ either vanishes everywhere or nowhere. We can proceed to the classification.

\section{Classification of Hypersurfaces.}
We will prove the following result first.
\begin{theorem}
Suppose $M\subset \Aff^{n+1}$ a hypersurface with a pointwise $SO(n-1)$ symmetry. Then its immersion $F$ is affine congruent to one of the following cases:
\begin{enumerate}
\item[1.] When $\zeta\neq 0$, then $F(t,\vec{u})$ is given by
\begin{equation}
F(t,\vec{u})=\left(\gamma_1(t)\phi(\vec{u}),\gamma_2(t)\right). \label{znon0}
\end{equation}
Here, $\phi$ is the immersion of an $\epsilon$-quadric in a hyperplane, where $\epsilon\neq 0$.
\item[2.] When $\zeta=0$ and $\sigma\neq r$, $F(t,\vec{u})$ is given by
\begin{equation}
F(t,\vec{u})=\left(\gamma_1(t)\vec{u},\gamma_1(t)\frac{\norm{\vec{u}}^2}{2}+\gamma_2(t),\gamma_1(t)\right). \label{rnots}
\end{equation}
\item[3.] When $\zeta=\sigma-r=0$, then $F(t,\vec{u})$ is given by
\begin{equation}
F(t,\vec{u})=\left(\vec{u},\frac{\norm{\vec{u}}^2}{2}+\gamma_2(t),\gamma_1(t)\right).
\label{riss}
\end{equation}
\end{enumerate}
\end{theorem}
When we've proven this, we will pick the immersions apart and check what conditions the planar curve $\gamma=(\gamma_1,\gamma_2)$ should satisfy such that the affine metric on the hypersurface is definite.
\begin{proof}
We start by noticing that for vector fields $X_i$ and $X_j$ orthogonal to $X_1$ we have:
\begin{equation*}
D_{X_i}X_j=\tilde{\nabla}_{X_i}X_j+\delta_{ij}\left(-(\sigma+r) X_1+\xi\right).
\end{equation*}
In this formula, we have $\tilde{\nabla}$ as the restriction of $\nabla$ to $\scal{X_1}^{\bot}$.
We consider the map $\phi=-(\sigma+r) X_1+\xi$ on $N$. By taking the derivative, we find
\begin{equation*}
\phi_*X_i=D_{X_i}\phi=-(\sigma+r)(\sigma-r)X_i-bX_i=-(\sigma^2-r^2+b)X_i=-\zeta X_i.
\end{equation*}
Calculating its derivative with respect to $X_1$, we find
\begin{align*}
D_{X_1}\phi&=-X_1(\sigma+r)X_1-(\sigma+r)\left((n-1)rX_1+\xi\right)-aX_1\\
           &=-(\sigma+r)\left(-(\sigma+r)X_1+\xi\right)=-(\sigma+r)\phi.
\end{align*}
It suffices to multiply $\phi$ with a factor $\beta\neq 0$ given as a solution to
\begin{equation}
\begin{split}
X_1(\beta)&=\beta (\sigma+r),\\
X_i(\beta)&=0,\label{beta}
\end{split}
\end{equation}
to make $\beta\phi$ independent of $X_1$. A straightforward calculation shows that the integrability conditions are satisfied, so a solution for $\beta$ exists. So either $\beta\phi$ is a constant vector along $M$, or $\beta\phi$ is an immersion of $N$ in $\Aff^{n+1}$. We rename $\beta \phi$ as simply $\phi$. \\\\
First we assume that $\zeta\neq 0$. In this case, we can regard $\phi$ as an immersion of $N$. We can calculate what structures it induces. We find
\begin{align*}
D_{X_i} \phi_*X_j&=-\beta \zeta D_{X_i}X_j\\
                 &=-\beta\zeta\left(\tilde{\nabla}_{X_i}X_j+\delta_{ij}\frac{\phi}{\beta}\right)\\
                 &=\phi_*(\tilde{\nabla}_{X_i}X_j)-\delta_{ij}\zeta \phi\\
                 &=\phi_*(\tilde{\nabla}_{X_i}X_j)-\epsilon h_N(X_i,X_j) \phi.
\end{align*}
This means that the immersion $\phi$ lies in a hyperplane of $\Aff^{n+1}$ spanned by $$\{X_2(p),\cdots,X_n(p),\phi(p)\}$$ for a starting point $p$, see \cite{erbacher}.  Hence, $\phi$ can be considered a hypersurface in $\Aff^n$. The position vector $\phi$ of this immersion is transversal to the tangent space, hence using it as a transversal allows the induced volume form to be parallel. We notice that the determinant $\det\left(h_N(X_i,X_j)\right)$ is a constant along the submanifold. We can then calculate
\begin{align*}
X_i\left(\omega_N(X_2,\cdots,X_n)\right)=\sum_{j=2}^n\omega_N(X_2,\cdots,\tilde{\nabla}_{X_i}X_j,\cdots,X_n)=0.
\end{align*}
Hence $\omega_N(X_2,\cdots,X_n)$ is also a constant, so $\phi$ is the affine normal up to a constant factor. We obtain that $\phi$ is an immersion of a proper affine hypersphere in a given hyperplane $H\subset\Aff^{n+1}$.
 Furthermore, $h_N$ is for every $t$ a constant multiple of $h$ and it is easy to see that the connection $\tilde{\nabla}$ coincides with its Levi-Civita connection. Hence $\phi$ is an immersion of $N$ as an $\epsilon$-quadric in $H$. This is already hinted by the fact that this hypersurface has a constant sectional curvature, see \cite{simonlizhao}. Next, we define a vector field along $N$ transversal to the hyperplane which contains the quadric. Therefore, we define $C$ as
\begin{equation*}
 C=bX_1+(\sigma-r)\xi.
\end{equation*}
We calculate its derivatives and find
\begin{align*}
D_{X_i}C&=b(\sigma-r)X_i-b(\sigma-r)X_i=0,\\
D_{X_1}C&=(nr-\sigma)C.
\end{align*}
Hence for every choice of $t$, $C$ is constant along $N$ and the direction of $C$ remains the same. We can see that $C$ is not contained within $H$, since we find that
\begin{equation*}
C=\lambda \frac{\phi}{\beta}+\sum_{j=2}^n a_j X_j \Leftrightarrow \zeta=0\bigwedge\forall j\in\{2,\cdots,n\}: a_j=0.
\end{equation*}
Since the derivative to $t$ of $C$ is parallel with $C$, there exists a non-vanishing function $c(t)$ and a constant vector $C_0$ such that $C(t)=c(t) C_0$. We can write $X_1$ as a linear combination of $C$ and $\phi$, thus eliminating $\xi$ and this looks like
\begin{equation*}
X_1=\pard{F}{t}=\frac{c}{\zeta}(t) C_0-\frac{\sigma-r}{\zeta\beta}(t)\phi.
\end{equation*}
Putting $C_0$ along $e_{n+1}$ and integrating we can define
\begin{align*}
\gamma_1(t)&=-\int\frac{\sigma-r}{\zeta\beta}(t) dt,\\
\gamma_2(t)&=\int \frac{c(t)}{\zeta(t)} dt.
\end{align*}
The resultant immersion, after applying an affine transformation is given by \eqref{znon0}.\\
Now we return to the case where $\zeta$ vanishes. Notice that in this case $X_1(\sigma-r)=(n+1)r(\sigma-r)$, so either $\sigma-r$ vanishes or $\sigma\neq r$ everywhere. In this case the defined vector field $\phi$ is a constant vector field on the hypersurface. We find
\begin{equation*}
D_{X_i}X_j=\tilde{\nabla}_{X_i}X_j+\delta_{ij} \frac{\phi}{\beta}.
\end{equation*}
Using the same method as when $\zeta\neq0$, we find that $N$ is immersed as a paraboloid in a hyperplane $H$. We can solve the equation
\begin{equation}
D_{X_1}X_1=(n-1)r X_1+\xi=(\sigma+nr)X_1+\frac{\phi}{\beta}.\label{diffeq2}
\end{equation}
Taking an appropriate initial condition for $\beta$, $\phi$ can be chosen as the constant affine normal $e_n$ of the paraboloid. We can solve \eqref{diffeq2} for $X_1$ and then integrate to obtain the immersion $F$. To find a unique solution, we can introduce initial conditions for $t=t_0$. For $F(t_0,\vec{u})$ we can use the paraboloid in the hyperplane $H\leftrightarrow x_{n+1}=0$, hence
\begin{equation}F(t_0,\vec{u})=\left(\vec{u},\frac{\norm{\vec{u}}^2}{2},0\right).\label{init1}\end{equation}
Then we consider $\left(X_1-(\sigma-r)F\right)(t_0)$, which is a constant, non-zero vector field along $N$.
This vector field is transversal to $H$. Applying an equi-affine transformation, we can place it along $e_{n+1}$. So there are numbers $\alpha_1$ and $\alpha_2$ for which $X_1(t_0)$ is given as
\begin{equation}
X_1(t_0)=\left(\alpha_1\vec{u},\alpha_1\frac{\norm{\vec{u}}^2}{2},\alpha_2\right),\label{init2}
\end{equation}
where $\alpha_1\alpha_2\neq0$.
The solution to \eqref{diffeq2} given \eqref{init1} and \eqref{init2} and applying an affine transformation is given by \eqref{rnots}.\\
Finally, when $\sigma-r$ vanishes, we can still consider the equation \eqref{diffeq2} with \eqref{init1}, but now $X_1(t_0)$ can be considered the constant transversal to $H$ and the initial condition becomes $X_1(t_0)=e_{n+1}$. The solution to \eqref{diffeq2} now is given by \eqref{riss} up to an affine transformation.
\end{proof}

For each case, we study what conditions we should impose on $\gamma$ such that the resultant hypersurface is a convex affine hypersurface. Note that the affine normal $\tilde{\phi}$ of the underlying quadric in the hyperplane is on itself transversal to this hypersurface. Denoting $\partial_t$ tangent to the part related to $\gamma$, a candidate for an affine normal is given as
\begin{equation}
\xi=\mu \tilde{\phi}+\nu \partial_t.\label{normal}
\end{equation}
We denote the coordinate vector fields tangent to the quadric by $\partial_{u_i}$.
\begin{theorem}
A hypersurface obtained by \eqref{znon0}, where $\phi$ is the immersion of an $\epsilon$-quadric ($\epsilon\neq 0$) is a convex affine hypersurface with a pointwise $SO(n-1)$ symmetry  when
\begin{equation}
\epsilon\dot{\gamma}_2\gamma_1(\ddot{\gamma}_1\dot{\gamma}_2-\ddot{\gamma}_2\dot{\gamma}_1)< 0.\label{metr1}
\end{equation}
Its affine normal is obtained through \eqref{normal} by
\begin{equation}
\begin{split}
\dot{\mu}&=-\nu\frac{\ddot{\gamma}_1\dot{\gamma}_2-\ddot{\gamma}_2\dot{\gamma}_1}{\dot{\gamma}_2},\\
\mu^{n+2}\gamma_1^{n-1}\dot{\gamma}_2^3&=\pm\epsilon^{n-1}\left(\ddot{\gamma}_1\dot{\gamma}_2-\ddot{\gamma}_2\dot{\gamma}_1\right).\label{sfeer1}
\end{split}
\end{equation}
\end{theorem}
\begin{proof}
We can calculate the derivatives of \eqref{znon0} and decompose them to find
\begin{align*}
\partial_t&=(\dot{\gamma}_1\phi,\dot{\gamma}_2),\\
\partial_{u_i}&=(\gamma_1\phi_{u_i},0),\\
D_{\partial_t}\partial_t&=\frac{\ddot{\gamma}_2}{\dot{\gamma}_2}\partial_t+\frac{\ddot{\gamma}_1\dot{\gamma}_2-\ddot{\gamma}_2\dot{\gamma}_1}{\dot{\gamma}_2}(\phi,0),\\
D_{\partial_t}\partial_{u_i}&=\frac{\dot{\gamma}_1}{\gamma_1} \partial_{u_i},\\
D_{\partial_{u_i}}\partial_{u_j}&=\tilde{\nabla}_{\partial_{u_i}}\partial_{u_j}-\epsilon g_{ij}\gamma_1 (\phi,0).
\end{align*}
Here $g_{ij}$ stands for the positive definite metric of $\phi$. For this to be a convex hypersurface, the induced affine metric must be definite and hence we find \eqref{metr1}. Calculating the derivatives of \eqref{normal} we find
\begin{align*}
D_{\partial_t}\xi&=\left(\dot{\nu}+\nu\frac{\ddot{\gamma}_2}{\dot{\gamma}_2}\right)\partial_t+\left(\dot{\mu}+\nu\frac{\ddot{\gamma}_1\dot{\gamma}_2-\ddot{\gamma}_2\dot{\gamma}_1}{\dot{\gamma}_2}\right)(\phi,0),\\
D_{\partial_{u_i}}\xi&=\frac{\mu+\nu\dot{\gamma}_1}{\gamma_1} \partial_{u_i}+\nu_{u_i} \partial_t+\mu_{u_i} (\phi,0).
\end{align*}
This transversal is equi-affine when the first equation in \eqref{sfeer1} is satisfied and satisfies the symmetry when $\mu$ and $\nu$ only depend on $t$. To find the affine normal, we rewrite the covariant derivatives and obtain
\begin{align*}
D_{\partial_t}\partial_t&=\left(\frac{\ddot{\gamma}_2}{\dot{\gamma}_2}-\frac{\nu}{\mu}\frac{\ddot{\gamma}_1\dot{\gamma}_2-\ddot{\gamma}_2\dot{\gamma}_1}{\dot{\gamma}_2}\right)\partial_t+\frac{1}{\mu}\frac{\ddot{\gamma}_1\dot{\gamma}_2-\ddot{\gamma}_2\dot{\gamma}_1}{\dot{\gamma}_2}\xi,\\
D_{\partial_t}\partial_{u_i}&=\frac{\dot{\gamma}_1}{\gamma_1} \partial_{u_i},\\
D_{\partial_{u_i}}\partial_{u_j}&=\tilde{\nabla}_{\partial_{u_i}}\partial_{u_j}+\frac{\nu}{\mu}\epsilon g_{ij} \partial_t-\epsilon g_{ij}\gamma_1 \frac{1}{\mu}\xi.
\end{align*}
The determinant of the affine metric is given as
\begin{equation*}
\det(h)=(-\epsilon\gamma_1)^{n-1}\frac{\ddot{\gamma}_1\dot{\gamma}_2-\ddot{\gamma}_2\dot{\gamma}_1}{\dot{\gamma}_2} \frac{1}{\mu^{n}}\det(g_{ij})
\end{equation*}
Calculating the volume form, we find
\begin{align*}
\omega(\partial_{u_2},\cdots,\partial_{u_n},\partial_t)&=\begin{vmatrix} \partial_{u_2} &\cdots &\partial_{u_n} & \partial_t &\mu(\phi,0)+\nu \partial_t \end{vmatrix}\\
                                                       &=-\begin{vmatrix} \gamma_1 \phi_{u_2} &\cdots &\gamma_1 \phi_{u_n}& \mu \phi &\dot{\gamma}_1 \phi\\ 0 & \cdots & 0 & 0 & \dot{\gamma}_2\end{vmatrix}\\
                                                       &=-\dot{\gamma}_2 \gamma_1^{n-1} \mu \begin{vmatrix} \phi_{u_2} & \cdots & \phi_{u_n} & \phi\end{vmatrix}.
\end{align*}
The condition for $\xi$ to be the affine normal is given by $$\abs{\det(h_{ij})}=\omega(\partial_{u_2},\cdots,\partial_{u_n},\partial_t)^2$$
which implies the final equation in \eqref{sfeer1}.\\
Finally, we check whether the cubic form satisfies the given symmetry and obtain
\begin{align*}
\nabla h(\partial_t;\partial_{u_i},\partial_{u_j})&=\pard{h(\partial_{u_i},\partial_{u_j})}{t}-h(\nabla_{\partial_t}\partial_{u_i},\partial_{u_j})-h(\nabla_{\partial_t}\partial_{u_i},\partial_{u_j})\\
                      &=-\left(\frac{\dot{\gamma}_1}{\gamma_1}+\frac{\dot{\mu}}{\mu}\right)h(\partial_{u_i},\partial_{u_j}),\\
\nabla h(\partial_t;\partial_{u_i},\partial_t)&=0,\\
\nabla h(\partial_{u_i};\partial_{u_j},\partial_{u_k})&=-\epsilon\gamma_1 \frac{1}{\mu}\nabla g(\partial_{u_i};\partial_{u_j},\partial_{u_k})=0 .
\end{align*}
This shows that the symmetry in $K$ is appropriate.\end{proof}
In case we look at the curve $\gamma=\left(e^{-\frac{t}{n}},e^{t}\right)$, the resultant hypersurface is the Calabi-type composition of a proper hyperbolic affine sphere and a point, see (\cite{dill}, section 3).
\begin{theorem}
A hypersurface obtained by \eqref{rnots} is a convex affine hypersurface with a pointwise $SO(n-1)$ symmetry  when
\begin{equation}
\dot{\gamma}_1\gamma_1(\ddot{\gamma}_2\dot{\gamma}_1-\ddot{\gamma}_1\dot{\gamma}_2)>0.\label{metr2}
\end{equation}
Its affine normal is obtained through \eqref{normal} by
\begin{equation}
\begin{split}
\dot{\mu}&=-\nu\frac{\ddot{\gamma}_2\dot{\gamma}_1-\ddot{\gamma}_1\dot{\gamma}_2}{\dot{\gamma}_1},\\
\mu^{n+2} \dot{\gamma}_1^3 \gamma_1^{n-1}&=\pm\left(\ddot{\gamma}_2\dot{\gamma}_1-\ddot{\gamma}_1\dot{\gamma}_2\right).\label{sfeer2}
\end{split}
\end{equation}
\end{theorem}
\begin{proof}
Taking the derivatives of \eqref{rnots} we obtain
\begin{align*}
\partial_t&=\left(\dot{\gamma}_1\vec{u},\dot{\gamma}_1\frac{\norm{\vec{u}}^2}{2}+\dot{\gamma}_2(t),\dot{\gamma}_1(t)\right),\\
\partial_{u_i}&=\gamma_1 e_{i-1}+u_i \gamma_1e_n,\\
D_{\partial_t} \partial_t&=\frac{\ddot{\gamma}_1}{\dot{\gamma}_1} \partial_t+\frac{\ddot{\gamma}_2\dot{\gamma}_1-\ddot{\gamma}_1\dot{\gamma}_2}{\dot{\gamma}_1} e_n,\\
D_{\partial_t}\partial_{u_i}&=\frac{\dot{\gamma}_1}{\gamma_1} \partial_{u_i},\\
D_{\partial_{u_i}}\partial_{u_j}&=\gamma_1 \delta_{ij} e_n.
\end{align*}
So in order to have a convex affine hypersurface, we get \eqref{metr2}.\\
Taking the derivatives of \eqref{normal} we find
\begin{align*}
D_{\partial_t}\xi&=\left(\dot{\nu}+\nu\frac{\ddot{\gamma}_1}{\dot{\gamma}_1}\right)\partial_t+\left(\dot{\mu}+\nu\frac{\ddot{\gamma}_2\dot{\gamma}_1-\ddot{\gamma}_1\dot{\gamma}_2}{\dot{\gamma}_1}\right) e_n,\\
D_{\partial_{u_i}}\xi&=\nu \frac{\dot{\gamma}_1}{\gamma_1} \partial_{u_i}+\pard{\nu}{u_i} \partial_t+\pard{\mu}{u_i} e_n.
\end{align*}
This is the desired affine normal when $\mu$ and $\nu$ only depend on $t$ and the first equation in \eqref{sfeer2} holds.
Rewriting the second derivatives of the immersion yields
\begin{align*}
D_{\partial_t}{\partial_t}&=\left(\frac{\ddot{\gamma}_1}{\dot{\gamma}_1}-\frac{\nu}{\mu}\frac{\ddot{\gamma}_2\dot{\gamma}_1-\ddot{\gamma}_1\dot{\gamma}_2}{\dot{\gamma}_1}\right) \partial_t+\frac{1}{\mu} \frac{\ddot{\gamma}_2\dot{\gamma}_1-\ddot{\gamma}_1\dot{\gamma}_2}{\dot{\gamma}_1} \xi,\\
D_{\partial_t}\partial_{u_i}&=\frac{\dot{\gamma}_1}{\gamma_1} \partial_{u_i},\\
D_{\partial_{u_i}}\partial_{u_j}&=\delta_{ij}\left(-\frac{\nu}{\mu} \gamma_1 \partial_t+\frac{1}{\mu}\gamma_1 \xi\right).
\end{align*}
The determinant of the affine metric in this basis is given by
\begin{equation*}
\det(h)=\frac{\left(\ddot{\gamma}_2\dot{\gamma}_1-\ddot{\gamma}_1\dot{\gamma}_2\right) \gamma_1^{n-1}}{\mu^n\dot{\gamma}_1}.
\end{equation*}
Calculating the volume form, we find
\begin{align*}
\omega(\partial_{u_2},\cdots,\partial_{u_n},\partial_t)&=\begin{vmatrix} \partial_{u_2} &\cdots & \partial_{u_n} & \partial_t & \mu e_n+\nu \partial_t \end{vmatrix}\\
                                  &=\begin{vmatrix} \gamma_1 \Eins_{n-1}  &\dot{\gamma}_1 \vec{u} & 0\\ \gamma_1 \vec{u}^t & \dot{\gamma}_1 \frac{\norm{\vec{u}}^2}{2}+\dot{\gamma}_2 & \mu \\ 0 & \dot{\gamma}_1 & 0 \end{vmatrix}\\
                                  &=-\mu \dot{\gamma}_1 \gamma_1^{n-1}.
\end{align*}
We can thus conclude that $\xi$ is an affine normal if the second equation in \eqref{sfeer2} holds.\\
Checking the affine metric, we find
\begin{align*}
\nabla h(\partial_t;\partial_{u_i},\partial_{u_j})&=-\left(\frac{\dot{\gamma}_1}{\gamma_1}+\frac{\dot{\mu}}{\mu}\right)h(\partial_{u_i},\partial_{u_j}),\\
\nabla h(\partial_t;\partial_{u_i},\partial_t)&=0,\\
\nabla h(\partial_{u_i};\partial_{u_j},\partial_{u_k})&=0 .
\end{align*}
\end{proof}
When we take $\gamma=\left(e^{\frac{t}{n+1}},e^{-\frac{t}{n+1}}\right)$, we obtain a semi-projective Calabi-type composition of an improper affine sphere and a point, see (\cite{dill},section 7).
\begin{theorem}
A hypersurface obtained by \eqref{riss} is a convex affine hypersurface with a pointwise $SO(n-1)$ symmetry  when
\begin{equation}
\left(\ddot{\gamma}_1\dot{\gamma}_2-\ddot{\gamma}_2\dot{\gamma}_1\right)\dot{\gamma}_2>0.\label{metr3}
\end{equation}
Its affine normal is obtained through \eqref{normal} by
\begin{equation}
\begin{split}
\dot{\mu}&=-\nu\frac{\ddot{\gamma}_1\dot{\gamma}_2-\ddot{\gamma}_2\dot{\gamma}_1}{\dot{\gamma}_2},\\
\mu^{n+2}\dot{\gamma}_2^3&=\pm\left(\ddot{\gamma}_1\dot{\gamma}_2-\ddot{\gamma}_2\dot{\gamma}_1\right).\label{sfeer3}
\end{split}
\end{equation}
\end{theorem}
\begin{proof}
We take the derivatives of \eqref{riss} and obtain
\begin{align*}
\partial_t&=\dot{\gamma}_1 e_n+\dot{\gamma}_2 e_{n+1},\\
\partial_{u_i}&=e_{i-1}+u_i e_n,\\
D_{\partial_t}\partial_t&=\frac{\ddot{\gamma}_2}{\dot{\gamma}_2} \partial_t+\frac{\ddot{\gamma}_1\dot{\gamma}_2-\ddot{\gamma}_2\dot{\gamma}_1}{\dot{\gamma}_2} e_n,\\
D_{\partial_t}\partial_{u_i}&=0,\\
D_{\partial_{u_i}}\partial_{u_j}&=\delta_{ij} e_n.
\end{align*}
For this to be a convex affine hypersurface we obtain \eqref{metr3}. We take the derivative of \eqref{normal} and obtain
\begin{align*}
D_{\partial_t}\xi&=\left(\dot{\nu}+\nu \frac{\ddot{\gamma}_2}{\dot{\gamma}_2}\right)\partial_t+\left(\dot{\mu}+\nu\frac{\ddot{\gamma}_1\dot{\gamma}_2-\ddot{\gamma}_2\dot{\gamma}_1}{\dot{\gamma}_2}\right)e_n,\\
D_{\partial_{u_i}}\xi&=\mu_{u_i} e_n+\nu_{u_i}\partial_t.
\end{align*}
The desired result is obtained when $\mu$ and $\nu$ only depend on $t$ and satisfy the first equation in \eqref{sfeer3}.\\
Rewriting the second derivatives we obtain
\begin{align*}
D_{\partial_t}\partial_t&=\left(\frac{\ddot{\gamma}_2}{\dot{\gamma}_2}-\frac{\nu}{\mu}\frac{\ddot{\gamma}_1\dot{\gamma}_2-\ddot{\gamma}_2\dot{\gamma}_1}{\dot{\gamma}_2}\right) \partial_t+\frac{1}{\mu}\frac{\ddot{\gamma}_1\dot{\gamma}_2-\ddot{\gamma}_2\dot{\gamma}_1}{\dot{\gamma}_2} \xi,\\
D_{\partial_t}\partial_{u_i}&=0,\\
D_{\partial_{u_i}}\partial_{u_j}&=\delta_{ij} \left(-\frac{\nu}{\mu} \partial_t+\frac{1}{\mu}\xi\right).
\end{align*}
The determinant of the affine metric is given as
\begin{equation*}
\det(h)=\frac{1}{\mu^n}\frac{\ddot{\gamma}_1\dot{\gamma}_2-\ddot{\gamma}_2\dot{\gamma}_1}{\dot{\gamma}_2}
\end{equation*}
and the volume form is given by
\begin{align*}
\omega(\partial_{u_2},\cdots,\partial_{u_n},\partial_t)&=\begin{vmatrix} \Eins_{n-1} & 0 & 0\\ \vec{u}^t & \dot{\gamma}_1 & \mu\\ 0& \dot{\gamma}_2 & 0 \end{vmatrix}\\
 &=-\mu\dot{\gamma}_2.
\end{align*}
Then $\xi$ is the affine normal when the final equation in \eqref{sfeer3} is satisfied.
Checking the affine metric, we find
\begin{align*}
\nabla h(\partial_t;\partial_{u_i},\partial_{u_j})&=-\frac{\dot{\mu}}{\mu}h(\partial_{u_i},\partial_{u_j}),\\
\nabla h(\partial_t;\partial_{u_i},\partial_t)&=0,\\
\nabla h(\partial_{u_i};\partial_{u_j},\partial_{u_k})&=0 .
\end{align*}
\end{proof}
When we take $\gamma=\left(e^{\frac{t}{2}},e^{-\frac{t}{n+1}}\right)$, then we obtain a standard Calabi-type composition of an improper affine sphere and a point, see (\cite{dill}, section 6).

\section{Regarding affine spheres.}
In the previous section, we obtained a classification for convex affine hypersurfaces and constructed their affine normals. Now we can specify what conditions are to be imposed in order to obtain affine spheres. We propose the general solution first and construct some examples later on.
\begin{theorem}
The improper affine spheres having a pointwise $SO(n-1)$ symmetry are affine congruent with one of the cases:
\begin{enumerate}
\item[1.] Immersion \eqref{znon0} where $\gamma$ satisfies
\begin{equation}c \dot{\gamma}_1^{n+2}\gamma_1^{n-1}=\dot{\gamma}_2^{n-1}(\ddot{\gamma}_1\dot{\gamma}_2-\ddot{\gamma}_2\dot{\gamma}_1),\label{improp1}
\end{equation}
for $c\in \Rl_0$.
\item[2.] Immersion \eqref{rnots} where $\gamma=(t,ct^{n+1})$ and $c \in \Rl_0$.
\item[3.] Immersion \eqref{riss} where $\gamma=(c t^n,t^{n+1})$ and $c\in \Rl_0$.
\end{enumerate}
\end{theorem}
\begin{proof}
For the first $2$ cases, the derivatives of \eqref{normal} have to vanish. Considering the first case, this implies
\begin{align*}
\frac{d(\nu\dot{\gamma}_2)}{dt}&=0,\\
\mu+\nu \dot{\gamma}_1&=0.
\end{align*}
Hence, we find that $\mu=-c_0\frac{\dot{\gamma}_1}{\dot{\gamma}_2}$. Substituting this into \eqref{sfeer1}, we get \eqref{improp1}.\\
In the second case we find that $\nu=0$ and hence $\mu=c_0$. Equation \eqref{sfeer2} can be simplified, reparametrising such that $\dot{\gamma}_1=1$. In that case we find
\begin{equation*}
\ddot{\gamma}_2=c_1t^{n-1},
\end{equation*}
for an appropriate constant $c_1\neq 0$. Integrating this twice and applying an affine transformation, we get that $\gamma=(t,ct^{n+1})$.\\
For the final case, we take a step back. In the construction of \eqref{riss}, we let $\xi=(0,\cdots,0,1)$. We can solve for $r=\sigma$ that
\begin{align*}
\dot{r}&=-(n+1)r^2\Rightarrow r=\frac{1}{(n+1)t},\\
\dot{\beta}&=2r\beta\Rightarrow \beta=c_0 t^{\frac{2}{n+1}}.
\end{align*}
Solving the equation
\begin{equation*}
\pard{F}{t}=\frac{\xi}{2r}-\frac{\phi}{2r\beta}=\left(0,\cdots,0,-\frac{n+1}{2c_0}t^{\frac{n-1}{n+1}},\frac{n+1}{2}t\right)
\end{equation*}
with initial condition
\begin{equation*}
F(t_0,\vec{u})=\left(\vec{u},\frac{\norm{\vec{u}}^2}{2},0\right)
\end{equation*}
and reparametrising yields the desired solution.
\end{proof}
To construct an example of \eqref{improp1}, we consider curves of the form $\left(e^{at},e^{bt}\right)$ where $a$ and $b$ are constants. Equation \eqref{metr1} implies that $\epsilon ab^2(a-b)<0$. Equation \eqref{improp1} implies that $\frac{ab^n (a-b)}{a^{n+2}}e^{\left(n(b-2a)\right)t}$ is a constant, hence $b=2a$. So we get an improper affine sphere, up to parametrisation by
\begin{equation*}
F(t,\vec{u})=\left(e^t \phi(\vec{u}),e^{2t}\right),
\end{equation*}
where $\phi$ represents an ellipsoid.\\
For proper affine spheres we obtain the following result.
\begin{theorem}
The proper affine spheres having a pointwise $SO(n-1)$ symmetry are affine congruent with one of the cases:
\begin{enumerate}
\item[1.] Immersion \eqref{znon0} where $\gamma$ is a centro-affine curve satisfying
\begin{equation}c\left(\dot{\gamma}_2\gamma_1-\dot{\gamma}_1\gamma_2\right)^{n+2}\gamma_1^{n-1}=\dot{\gamma}_2^{n-1}\left(\ddot{\gamma}_1\dot{\gamma}_2-\ddot{\gamma}_2\dot{\gamma}_1\right),\label{prop1}
\end{equation}
for $c\in \Rl_0$.
\item[2.] Immersion \eqref{rnots} where $\gamma$ is a centro-affine curve satisfying
\begin{equation}c \left(\dot{\gamma}_1\gamma_2-\dot{\gamma}_2\gamma_1\right)^{n+2} \gamma_1^{n-1}=\dot{\gamma}_1^{n-1}\left(\ddot{\gamma}_2\dot{\gamma}_1-\ddot{\gamma}_1\dot{\gamma}_2\right),\label{prop2}
\end{equation}
for $c\in \Rl_0$.
\end{enumerate}
\end{theorem}
\begin{proof}
To prove this, we simply shift the affine normal to the position vector and \eqref{normal} gives the required $\mu$. Substituting this into \eqref{sfeer1} and \eqref{sfeer2} yields the result.
\end{proof}
We can construct examples of these cases using $\gamma=(e^{at},e^{bt})$. In the first case substitution in \eqref{prop1} yields
\begin{equation*}
c=-\frac{a b^n}{(b-a)^{n+1}} e^{-2(na+b)t}.
\end{equation*}
For the metric, we get $\epsilon ab^2(a-b)=\epsilon a^4 n^2(1+n)<0$. So $\phi$ is the immersion of a hyperboloid. The result is a hyperbolic affine sphere. In particular, it's the Calabi product of a point and a hyperboloid, see \cite{dill}. It's given as
\begin{equation*}
F(t,\vec{u})=\left(e^{t}\phi(\vec{u}),e^{-nt}\right).
\end{equation*}
To obtain an example of the final case, we use \eqref{prop2} on the considered curve which yields
\begin{equation*}
c=-\frac{a^nb}{(a-b)^{n+1}} e^{-(n+1)(a+b)t}.
\end{equation*}
Equation \eqref{metr2} requires that $2a^4>0$, which is satisfied for every choice of $a\neq 0$. The immersion
\begin{equation*}
F(t,\vec{u})=\left(e^t\vec{u},e^t\frac{\norm{\vec{u}}^2}{2}+e^{-t},e^t\right)
\end{equation*}
is a convex hyperbolic affine sphere.

 \end{document}